\documentclass[final,5p,times,twocolumn,authoryear]{elsarticle}

\usepackage{amssymb}
\usepackage{amsfonts}
\usepackage{latexsym, amssymb,amsmath, amsbsy,amsopn, amstext}
\usepackage{amsmath, amsbsy}
\usepackage{hyperref}
\usepackage{float}
\usepackage{ifpdf}
\usepackage{xcolor}
\usepackage{tikz}
\usetikzlibrary{arrows,shapes,snakes,shadows,positioning,automata,patterns}
\usetikzlibrary{trees,decorations.pathmorphing,decorations.markings}

\def\0{{\bf 0}}

\def\thefootnote{*}

\newtheorem{thm}{Theorem}[section]
\newtheorem{lem}[thm]{Lemma}

\newtheorem{exa}[thm]{Example}
\newtheorem{rem}[thm]{Remark}

\newtheorem{assum}{Assumption}

\journal{Automatica}
\begin{document}
\begin{frontmatter}

\title{Initialization-free Distributed Algorithms for Optimal Resource Allocation with Feasibility Constraints and Application to Economic Dispatch of Power Systems\tnoteref{label0}}

\author[label1]{Peng Yi}
\address[label1]{Department of Electrical \& Computer Engineering, University of Toronto}
\ead{peng.yi@utoronto.ca}

\author[label2]{Yiguang Hong \corref{cor1}}
\address[label2]{Key Lab of Systems and Control, Academy of
Mathematics and Systems Science, Chinese Academy of Sciences,}
\ead{yghong@iss.ac.cn}

\author[label3]{Feng Liu}
\address[label3]{Department of Electrical Engineering, Tsinghua University }
\ead{lfeng@mail.tsinghua.edu.cn}

\begin{abstract}
In this paper, the distributed resource allocation optimization
problem is investigated.
The allocation decisions are made to minimize the
sum of all the agents' local objective functions while satisfying
both the global network resource constraint and the local allocation
feasibility constraints. Here the data corresponding to each agent
in this separable optimization problem, such as the network
resources, the local allocation feasibility constraint, and the
local objective function, is only accessible to individual agent and
cannot be shared  with others, which renders new challenges in this
distributed optimization problem.   Based on either projection or
differentiated projection, two classes of continuous-time
algorithms are proposed to solve this distributed optimization
problem in an initialization-free and scalable manner. Thus, no
re-initialization is required even if the operation environment or
network configuration is changed, making it possible to achieve a
``plug-and-play" optimal operation of networked heterogeneous
agents.  The algorithm convergence is guaranteed for strictly convex
objective functions, and the exponential convergence is proved for
strongly convex functions without local constraints.
Then the proposed algorithm is applied
to the distributed economic dispatch problem in power grids, to
demonstrate how it can achieve the global optimum in a scalable
way, even when the generation cost, or system load, or network
configuration, is changing.
\end{abstract}

\tnotetext[label0]{This paper was not presented at any IFAC
meeting. This work is supported by Beijing Natural Science
Foundation (4152057), NSFC (61333001, 61573344), and Program 973
(2014CB845301/2/3).
This work is also partly supported by the National Natural Science Foundation
of China (No. 51377092), Foundation of Chinese Scholarship Council (CSC No. 201506215034).
Corresponding author: Yiguang Hong,  Tel. +86(010)82541824.
Fax +86(010)62587343.}
\end{frontmatter}



\section{Introduction}

Resource allocation is one of the most important problems in network
optimization, which has been widely investigated in various areas
such as economics systems, communication networks, sensor networks,
and power grids. The allocation decisions may be made centrally by
gathering all the network data together to a decision-making center,
and then sent back to corresponding agents (referring to
\cite{RA0}). On the other hand, differing from this centralized
policy, the master-slave-type decentralized  algorithms, either
price-based (\cite{RA1}) or resource-based (\cite{RA2}), are
constructed to achieve the optimal allocations by the local
computations in the slave agents under the coordinations of the
master/center  through a ¡°one-to-all¡± communication architecture.
However, these methods may not be suitable or effective for the
resource allocation in large-scale networks with numerous
heterogeneous agents due to complicated network structures, heavy
communication burden, privacy concerns, unbearable time delays, and
unexpected single-point failures. Therefore, fully distributed
resource allocation optimization algorithms are highly desirable.

Distributed optimization, which cooperatively achieves optimal
decisions by the local manipulation with private data  and the
diffusion of  local information through a multi-agent network, has
drawn more and more research attention in recent years.  To
circumvent the requirement of control center or master, various
distributed optimization models or algorithms have been developed
(\cite{Ned2}, \cite{sayd}, \cite{lou12}, and
\cite{peng2}). In light of the increasing attention to distributed
optimization  and the seminal  work on distributed  resource
allocation in  \cite{RA3}, some distributed algorithms for resource
allocation optimization have been proposed  in \cite{RA4},
\cite{RA5}, \cite{RA6}, \cite{RA7}, and \cite{RA8}.

Continuous-time gradient flow algorithms have been widely
investigated for convex optimization after the pioneer work
\cite{arrow},  and  detailed references can be found in
\cite{arrow2} and \cite{bhaya}. Gradient flow algorithms  have been
applied to  network control and optimization ( \cite{CA1},
\cite{CA2} and \cite{CA3}), neural networks (\cite{arrow2}), and
stochastic approximation (\cite{proddup2}).  Recently,
continuous-time gradient flow algorithms have been adopted for
solving unconstrained distributed optimization problems (see
\cite{wang}, \cite{cort2}, \cite{mug2}, and \cite{cor1}).
Furthermore,  the projection-based gradient flow dynamics have
been employed for solving the complicated constrained optimization problems
in \cite{aubin2},  \cite{prods2}, \cite{jun}, \cite{Gao} and
\cite{cortes_prods}, and  the projected gradient flow ideas began to be  applied to distributed constrained optimization (see
\cite{liu}, \cite{xie} and \cite{peng3}).

 The economic dispatch, one of the key concerns
in power grids,  is to find the optimal secure
generation allocation to balance the system loads, and hence, can be regarded as
a special resource allocation problem. In recent years, there has been increasing research
attention in solving economic dispatch problems through a multi-agent
system in a distributed manner to meet the ever growing challenges
raised by increasing penetration of renewable energies and deregulation of power infrastructure (\cite{ED_zam} and \cite{ED2}).
 Mathematically, this boils down to a particular  distributed resource allocation
optimization problem. Furthermore, there were various continuous-time algorithms for the Distributed Economic Dispatch Problem (DEDP). For example,
\cite{ED1} showed that the physical power grid dynamics could serve as a part of a primal-dual gradient flow
algorithm to solve the DEDP, and in fact, it
considered physical network interconnections and generator dynamics explicitly,
 providing a quite comprehensive method and inspiring insights.   Moreover,
\cite{ED3} solved the DEDP  by combining the penalty method and the
distributed continuous-time algorithm in \cite{RA3}, and  proposed a procedure to
fulfill the initialization requirement, while \cite{EDcc} constructed a novel initialization-free distributed
algorithm to achieve DEDP given one
agent knowing the total system loads.

Motivated by various practical problems, including the DEDP in power grids, we study a
Distributed Resource Allocation Optimization (DRAO) problem, where
each agent can only manipulate its private data, such as the
local objective function, Local Feasibility Constraint
(LFC), and local resource data. Such data in practice cannot be shared or known by other
agents.
As the total network resource is the sum of individual agent's local
resources, the agents need to cooperatively achieve the optimal
resource allocation in a distributed way, so that the global
objective function (as the sum of all local objective functions) is
minimized with all the  constraints (including the network resource
constraint and LFCs) satisfied.
Note that the LFC is critical for the (secure) operation of practical
networks (referring to the communication system in \cite{RA10} and
\cite{RA9} as an example), even though it was not considered in most existing DRAO
works. Particularly, for the DEDP in power grids,
the generation of each generator must be limited within its box-like
capacity bounds.
The consideration of LFCs brings remarkable difficulties to existing
distributed algorithms designed for the DRAO without
LFCs, because the KKT (optimality) conditions for the DRAO with and
without LFCs are totally different (referring to Remark
\ref{kktcompare}).
So far, many DEDP works (such as \cite{ED1}, \cite{ED2} and \cite{ED3} and \cite{EDcc}) have only considered the box-like LFCs. However, the requirement from power industries, such as the secure operation of inverter-based devices in smart grids, promotes the demand to deal with non-box LFCs.
This extension is nontrivial, and we will show how to handle it systematically by using projected dynamics in this paper.

Another crucial albeit difficult problem is the initialization
coordination among all agents. Many existing results
are based on initialization coordination procedures to guarantee
that the initial allocations satisfy the network resource
constraint, which may only work well for ¡°static¡± networks.
However, for a ¡°dynamical¡± network, the resource has to be
re-allocated once the network configuration changes.
Therefore, the initialization coordination has to be
re-performed whenever these optimization algorithms re-start, which
considerably degrades their applicability. Taking the DEDP as an
example, the initialization needs to be coordinated among all agents
whenever local load demand or generation capacity/cost changes, or any
distributed generator plugs in or leaves off (see
\cite{ED3} for an initialization procedure).
This issue has to be well addressed for achieving highly-flexible power grids with the integration
of ever-increasing renewables.

The objective of this paper is to propose  an initialization-free
methodology to solve the DARO with local LFCs.    The main technical
contributions of this paper are highlighted as follows:
\begin{itemize}
\item  By employing the (differentiated) projection operation,  two fully distributed  continuous-time algorithms are proposed as a kind of projected dynamics, with the local allocation of each agent kept within its own LFC set. Moreover, the algorithms ensure the network resource constraint asymptotically without requiring it being satisfied at the initial points.
      Therefore,  it is  initialization-free,  different  from those given in  \cite{RA4}, \cite{RA5}, \cite{RA6} and \cite{RA7},
      and moreover, provides novel initialization-free algorithms different from the one  given in  \cite{EDcc}.
  \item The  convergence of the two projected algorithms is shown by the properties
       of Laplacian matrix and projection operation as well as the LaSalle invariance principle.
       The result can be regarded as an extension of some existing distributed optimization algorithms ( \cite{wang}, \cite{cor1}, \cite{xie}, and \cite{liu}) and an application of projected dynamics for variational inequalities (\cite{Gao} and \cite{jun})  to the DRAO problem.
  \item The proposed algorithms can be directly applied to the DEDP in power systems  considering generation capacity limitations.
      It enables  the ¡°plug-and-play¡± operation for power grids with high-penetration of flexible renewables. Our algorithms are essentially different from the ones provided in \cite{ED3} and \cite{EDcc},  and address multi-dimensional decision variables and  general non-box LFCs.
      Simulation results demonstrate that the algorithm  effectively deals with various data and network configuration changes, and also illustrate the algorithm scalability.
  \end{itemize}

The reminder of this paper is organized as follows.
The preliminaries are given and then the DRAO with LFCs is formulated
with the basic  assumptions in Section 2. Then a distributed
algorithm  in the form of  projected dynamics is proposed with its
convergence  analysis in Section 3. In Section 4, a differentiated projected algorithm is proposed with its convergence  analysis for DRAO with strongly convex objective functions, and an  exponential convergence rate is obtained in the case without LFCs. Moreover, the
application to the DEDP in power systems is shown in section 5 with
numerical experiments. Finally, the concluding remarks are given in
Section 6.

Notations: Denote $\mathbf{R}_{\geq 0}$ as the set of nonnegative real numbers.
Denote $\mathbf{1}_m=(1,...,1)^T \in \mathbf{R}^m$ and
$\mathbf{0}_m=(0,...,0)^T \in \mathbf{R}^m$.
Denote $col(x_1,....,x_n)=(x_1^T,\cdots, x_n^T)^T $ as the column vector stacked with vectors $x_1,...,x_n$. $I_n$ denotes the identity matrix in $\mathbf{R}^{n\times n}$.
For a matrix $A=[a_{ij}]$, $a_{ij}$ or $A_{ij}$
stands for the matrix entry in the $i$th row and $j$th column of $A$.
$A \otimes B$ denotes the  Kronecker  product of matrixes $A$ and $B$.
Denote $\times_{i=1,...,n}\Omega_i$ as the Cartesian product of the sets $\Omega_i,i=1,...,n$.
Denote the set of interiors of set $K$ as $int(K)$, and the boundary set of set $K$ as $\partial K$.

\section{Preliminaries and problem formulation}

In this section, we first give the preliminary knowledge related to convex analysis and graph theory, and then formulate the DRAO problem of interest.

\subsection{Convex analysis and projection}

The following concepts and properties about convex functions, convex sets, and projection operations come from \cite{ob} and \cite{bersekas}.
A differentiable  convex function $f: \mathbf{R}^m\rightarrow \mathbf{R}$
has the locally Lipschitz continuous gradient, if, given any compact set $Q$,
there is a constant $k^Q$ such that
$ ||\nabla f(x) - \nabla f(y) ||   \leq   k^Q ||x-y ||, \forall x, y \in Q. $
A  differentiable  function $f(x)$ is called $\mu$-strongly convex on $\mathbf{R}^m$  if there
exists a constant $\mu >0$ such that, for any $x,y \in \mathbf{R}^m$,
$
f(y)\geq f(x)+ \nabla^T f(x)(y-x) + \frac{1}{2}\mu ||y-x ||^2,
$
or equivalently,
\begin{equation}\label{stronglyconvex}
(x-y)^T (\nabla f(x)- \nabla f(y)) \geq \mu ||x-y ||^2.
\end{equation}

The following notations describe the geometry properties of the convex  set.
Denote $C_{\Omega}(x)$ as the normal cone of $\Omega$ at $x$,  that is,
$C_{\Omega}(x) =\{v: \langle v, y-x\rangle \leq 0, \quad \forall y\in \Omega\}. $
Define $c_{\Omega}(x)$ as
$c_{\Omega}(x)=\{v: ||v||=1, \langle v, y-x \rangle \leq 0, \; \forall y\in \Omega \}$ if $x\in \partial \Omega$, and $c_{\Omega}=\{\mathbf{0}\}$ if $x\in int(\Omega)$.
The feasible direction cone of $\Omega$ at  $x$ is given as
$K_{\Omega}(x)= \{d: d=\beta (y-x),\;y\in \Omega,\; \beta\geq 0\}.  $
The tangent cone of set $\Omega$ at  $x$ is defined as
$T_{\Omega}(x)=\{ v: v = \lim_{k\rightarrow \infty} \frac{x^k-x}{\tau_k}, \tau_k\geq0, \tau_k\rightarrow 0, x^k\in \Omega, x^k\rightarrow x \}.$
Then $T_{\Omega}(x)$ is the closure of $K_{\Omega}(x)$ when $\Omega$ is a  closed convex set,
and is the polar cone to $C_{\Omega}(x)$, that is, $T_{\Omega}(x)=\{y: \langle y, d\rangle \leq 0, \forall d\in C_{\Omega}(x) \}$ (referring to Lemma 3.13 of \cite{ob}).

Define the projection of $x$ onto a
closed convex set $\Omega$ by $P_{\Omega}(x)=\arg\min_{y\in \Omega} ||x-y ||$.
 The basic property of projection operation is
\begin{equation}\label{projection}
\langle x-P_{\Omega}(x), P_{\Omega}(x)-y \rangle \geq 0, \forall x\in \mathbf{R}^m,  \forall y\in \Omega.
\end{equation}
The following relationships can be derived from \eqref{projection},
\begin{equation}\label{projection2}
||x-P_{\Omega}(x) ||_2^2 + ||P_{\Omega}(x)-y||_2^2 \leq || x-y||_2^2,  \forall x\in \mathbf{R}^m,  \forall y\in \Omega,
\end{equation}
and
\begin{equation}\label{projection3}
||P_{\Omega}(x)-P_{\Omega}(y) || \leq ||x-y ||, \forall x, y \in \mathbf{R}^m.
\end{equation}
The normal cone $C_{\Omega}(x)$ can also be defined as (Lemma 2.38 of \cite{ob})
\begin{equation}
C_{\Omega}(x)=\{v: P_{\Omega}(x+v)=x\}.
\end{equation}

For a closed convex set $\Omega$, point $x\in \Omega$ and direction $v$, we define the differentiated  projection operator as (\cite{proddup2} and \cite{prods2}),
\begin{equation}\label{dpoperator}
\Pi_{\Omega}(x, v) = \lim_{\delta\rightarrow 0}\frac{P_{\Omega}(x+\delta v)-x}{\delta}.
\end{equation}
The basic properties of the differentiated projection operator are given as  follows (\cite{prods4}).

\begin{lem}\label{dpoperatorbaic}
(i):If $x \in int(\Omega)$, then $\Pi_{\Omega}(x,v)=v$;
(ii): $x \in \partial \Omega$, and $\max_{n \in c_{\Omega}(x)} \langle v, n \rangle \leq 0 $, then $\Pi_{\Omega}(x,v)=v$;
(iii): $x \in \partial \Omega$, and $\max_{n \in c_{\Omega}(x)} \langle v, n \rangle \geq 0 $, then
     $\Pi_{\Omega}(x,v)=v-\langle v, n^*\rangle n^*,$ where
     $n^*=\arg\max_{n\in c_{\Omega}(x)} \langle v, n\rangle$.
Therefore, the operator $\Pi_{\Omega}(x,v)$ in \eqref{dpoperator} is equivalent with the projection of $v$ onto $T_{\Omega}(x)$, i.e.,
$$ \Pi_{\Omega}(x,v)=P_{T_{\Omega}(x)}(v).$$
\end{lem}

\subsection{Graph theory}

The following concepts of graph theory can be found in \cite{god}.
The information sharing or exchanging among the agents is described by  graph
$\mathcal{G}=(\mathcal{N},\mathcal{E})$.
The edge set $\mathcal{E} \subset \mathcal{N}\times \mathcal{N} $ contains all the  information interactions.
If agent $i$ can get  information from agent $j$, then $(j,i) \in \mathcal{E}$ and
agent $j$ belongs to agent $i$'s  neighbor set $\mathcal{N}_i=\{ j | (j,i) \in
\mathcal{E}\}$.
$\mathcal{G}$ is said to be undirected when
$(i,j)\in \mathcal{E}$ if and only if $(j,i)\in \mathcal{E}$.
A path of graph $\mathcal{G}$ is a sequence of distinct agents in
$\mathcal{N}$ such that any consecutive agents in the sequence
corresponding to an edge of graph $\mathcal{G}$. Agent $j$ is said to be
connected to agent $i$ if there is a path from $j$ to $i$.
$\mathcal{G}$ is said to be  connected if any two agents are
connected.

Define adjacency matrix $A=[a_{ij}]$ of $\mathcal{G}$ with $a_{ij}=1$ if
$j\in \mathcal{N}_i$ and $a_{ij}=0$ otherwise.
Define the degree matrix
$Deg=diag\{ \sum_{j=1}^n a_{1j},..., \sum_{j=1}^n a_{nj}\}.$
Then the Laplacian of graph $\mathcal{G}$ is
$L=Deg-A.$
When $\mathcal{G}$ is a  connected  undirected graph, 0 is a simple
eigenvalue of Laplacian $L$ with the eigenspace $\{ \alpha \mathbf{1}_n| \alpha\in \mathbf{R}\}$,  and
$L \mathbf{1}_n=\mathbf{0}_n$, $\mathbf{1}^T_{n} L=\mathbf{0}^T_n$, while all other eigenvalues are positive.
Denote the eigenvalues of $L$   in an ascending order as $0<s_2\leq \cdots \leq s_n$.  Then, by the Courant-Fischer Theorem,
\begin{equation}\label{lap}
\min_{x\neq \mathbf{0},\atop{\mathbf{1}^Tx=0}} x^T L x =s_2||x||_2^2, \quad \max_{x\neq \mathbf{0}} x^T L x = s_n||x ||^2_2.
\end{equation}

\subsection{Problem formulation}\label{sec1}

Consider a group of agents with the index set $\mathcal{N}=\{1,...,n\}$ to
make an optimal allocation of network resource
under both the network resource constraint and LFCs.
Agent $i$ can decide its local allocation $x_i \in \mathbf{R}^m$,
 and  can access the local resource data $d_i \in \mathbf{R}^m$.
The total  network resource  is $\sum_{i\in \mathcal{N}} d_i$,
and therefore, the allocation should satisfy the {\bf network resource constraint:
 $\sum_{i\in \mathcal{N}} x_i= \sum_{i\in \mathcal{N}} d_i$.}
Furthermore, the  allocation  of agent $i$ should satisfy the {\bf local  feasibility constraint (LFC): $x_i\in \Omega_i$,}
where $\Omega_i \subset \mathbf{R}^m$ is a closed convex set only known by agent $i$.
Agent $i$ also has a
local  objective function $f_i(x_i): \mathbf{R}^m \rightarrow \mathbf{R} $ associated with its local allocation $x_i$.
Denote $X=col(x_1,...,x_n)\in \mathbf{R}^{mn}$ as the  allocation vector of the whole network.
Then the task for the agents  is  to collectively find the optimal  allocation corresponding to the
DRAO problem  as follows:
\begin{equation}\label{CRA1}
\begin{array}{l}\hline
{\bf Distributed \; Resource  \; Allocation \; Optimization }:           \\ \hline
\displaystyle \min_{x_i\in \mathbf{R}^m, \; i\in \mathcal{N}}             f(X)     =\sum_{i\in \mathcal{N}} f_i(x_i)         \\
\displaystyle subject \; to \;                            \sum_{i\in \mathcal{N}} x_i=\sum_{i\in\mathcal{ N}} d_i,  \\
\displaystyle  \qquad \qquad               \;  x_i\in \Omega_i, \;  i\in \mathcal{N}.                              \\  \hline
\end{array}
\end{equation}

Clearly, problem \eqref{CRA1} is an extension of  the previous optimization models in \cite{RA4}, \cite{RA5} and \cite{RA6} by introducing
the additional LFCs, that is, $x_i \in \Omega_i$.  Clearly, $x_i\in \Omega_i$ also generalizes previous box constraints in \cite{RA10}, \cite{RA9} and \cite{EDcc}.

The following assumptions are given for (\ref{CRA1}), which were also adopted for the distributed optimization or resource allocation in \cite{CA2}, \cite{cor1}, and \cite{liu}.

\begin{assum}\label{asumFun}
The functions $f_i(x_i),i\in \mathcal{N}$ are continuously differentiable  convex functions with locally Lipschitz continuous gradients and  positive definite Hessians over $\mathbf{R}^m$.
\end{assum}

Assumption \ref{asumFun} implies that $f_i(x_i)$'s are strictly convex, and hence guarantees the uniqueness of the optimal solution to (\ref{CRA1}).

\begin{assum}\label{slater}
There exists a finite optimal
solution $X^*$ to  problem \eqref{CRA1}. The Slater's constraint condition is satisfied for DRAO (\ref{CRA1}), namely,
there exists $\tilde{x}_i \in int(\Omega_i),\forall i\in \mathcal{N}$,
such that  $\sum_{i\in \mathcal{N}} \tilde{x}_i=\sum_{i\in\mathcal{ N}} d_i. $
\end{assum}

\begin{rem}
Define the recession cone of a  convex set $\Omega$ as
$R_{\Omega}=\{d: x+\alpha d\in \Omega,  \; \forall \; \alpha \geq 0, \; \forall x\in \Omega \}. $
Then the sufficient and necessary condition for the existence of finite optimal solution to \eqref{CRA1} is  (referring to Proposition 3.2.2 of \cite{bersekas})
\begin{equation}\label{finiteness}
\times_{i\in \mathcal{N}}  R_{\Omega_i}
\cap \{ Null(\mathbf{1}^T_n\otimes I_m)\}
\cap \times_{i\in \mathcal{N}} R_{f_i}=\mathbf{0}, \nonumber
\end{equation}
where $R_{\Omega_i}$ is  the recession cone of $\Omega_i$ and $R_{f_i}$
 is the recession cone of  any  nonempty level set of $f_i(x_i)$: $\{x_i\in \mathbf{R}^m | f_i(x_i)\leq \gamma\}$.

In many practical cases,  we have $R_{f_i}=\mathbf{0}$ (taking the quadratic function as an example).
Furthermore,  $R_{\Omega_i}=\mathbf{0}$ when $\Omega_i$ is compact.
Therefore,  the existence of a finite solution can be easily guaranteed and verified in many practical problems.
\end{rem}

The local objective function $f_i(x_i)$, resource data $d_i$  and  LFC set $\Omega_i$
are   the {\bf private data} for agent $i$, which are not shared with other agents.
This makes (\ref{CRA1}) a distributed optimization  problem.
To fulfill the  cooperations between agents for solving \eqref{CRA1},
the agents have to share their local information through a  network $\mathcal{G}=(\mathcal{N},\mathcal{E})$.
Next follows an assumption  about the connectivity of $\mathcal{G}$ to guarantee that any agent's
information can reach any other agents, which is also quite standard for distributed optimization (\cite{liu}).

\begin{assum}\label{assum2}
The  information sharing  graph $\mathcal{G}=(\mathcal{N},\mathcal{E})$ is  undirected and  connected.
\end{assum}

In a sum, the task of this paper is to design { fully distributed  algorithms } for the agents
to cooperatively find  the optimal resource allocation  to (\ref{CRA1}) without any center.
In other words, agent $i$ needs to find its optimal  allocation $x_i^*$
by  manipulating  its local private data $d_i$, $\Omega_i$, and  $f_i(x_i)$
and  by  cooperations  with its neighbor agents through $\mathcal{G}$.

\section{Projected algorithm for DRAO}

In this section,  a distributed algorithm for \eqref{CRA1} based on  projected dynamics is
proposed and analyzed.
The distributed algorithm for agent $i$ is given as follows:
\begin{equation}\label{gcd}
\begin{array}{l}\hline
{\bf Projected \;  algorithm \; for \; agent \; i}:           \\ \hline
\displaystyle \dot{x_i}       = P_{\Omega_i}( x_i-  \nabla f_i(x_i) +\lambda_i)-x_i            \\
\displaystyle \dot{\lambda}_i =-\sum_{j\in \mathcal{N}_i}(\lambda_i-\lambda_j)-\sum_{j\in \mathcal{N}_i}(z_i-z_j)+(d_i-x_i)  \\
\displaystyle \dot{z_i}       =\sum_{j\in \mathcal{N}_i}(\lambda_i-\lambda_j)                              \\  \hline
\end{array}
\end{equation}

 In algorithm \eqref{gcd}, $x_i\in \mathbf{R}^m$ is the local allocation of agent $i$, and $\lambda_i,z_i \in \mathbf{R}^m$ are two auxiliary variables of agent $i$.
Note that the algorithm (\ref{gcd}) is  fully distributed because agent  $i$ only needs the  local data (including $d_i$, $f_i(x_i)$, and  $\Omega_i$)
and the shared  information $\{\lambda_j,z_j, j\in \mathcal{N}_i\}$  from its neighbor agents.
Thereby, \eqref{gcd} does not need any center to handle all the data or coordinate all the agents.
With the distributed algorithm \eqref{gcd}, each agent has the autonomy and authority to formulate its
own objective function and feasibility set, and hence, the privacy is kept within each agent.
Because each agent can instantaneously react to its local data changes, it can quickly adapt its local decision.
Therefore,  the algorithm can be  easily applied  to large-scale networks.

The algorithm (\ref{gcd}) can be  understood based on the following observations.
The duality of (\ref{CRA1}) with multiplier $\lambda \in \mathbf{R}^m$ is
\begin{equation}\label{problem2}
\max_{\lambda \in \mathbf{R}^m} q(\lambda)=\sum_{i\in \mathcal{N}} q_i(\lambda)= \sum_{i\in \mathcal{N}}\inf_{x_i\in \Omega_i}\{f_i(x_i)-\lambda^T x_i+\lambda^T d_i\}.
\end{equation}
Although some existing distributed algorithms in \cite{Ned2} and \cite{sayd} addressed the dual problem (\ref{problem2}),
they need to solve a subproblem at each time (iteration) to calculate the  gradients.
In other words, two ``time scales" are  needed if applying existing distributed algorithms to \eqref{CRA1}.
Here we aim to develop a simple algorithm  without solving any subproblems.
To this end, we formulate a constrained optimization  problem with Laplacian matrix $L$ and $\Lambda=col(\lambda_1,...,\lambda_n)\in \mathbf{R}^{mn}$ as
\begin{equation}\label{problem3}
\begin{array}{lll}
&\max_{\Lambda=(\lambda_1,...,\lambda_n)} \quad \;& Q(\Lambda)\; = \; \sum_{i\in \mathcal{N}} q_i(\lambda_i)\\
& subject \; to \; \quad\;& (L \otimes I_m) \Lambda=\mathbf{0}_{mn}.
\end{array}
\end{equation}
The augmented Lagrangian duality of (\ref{problem3}) with  multipliers $Z=col(z_1,...,z_n) \in \mathbf{R}^{mn}$ is
\begin{equation}\label{problem4}
\min_{Z} \max_{\Lambda} Q(\Lambda,Z)=\sum_{i\in \mathcal{N}} q_i(\lambda_i)  - Z^T (L\otimes I_m) \Lambda - \frac{1}{2} \Lambda^T (L\otimes I_m) \Lambda.
\end{equation}
Then the projected dynamics (\ref{gcd}) is derived by applying the gradient flow  to (\ref{problem2})  and (\ref{problem4}) along with the projection operation to guarantee the feasibility of LFCs.

From the KKT condition, we can show that the equilibrium point of \eqref{gcd} yields the optimal solution to  problem \eqref{CRA1}.
Denote $X=col(x_1,...,x_n)$, $\nabla f(X)=col(\nabla f_1(x_1),\cdots,\nabla f_n(x_n))$, $D=col(d_1,\cdots,d_n)$,
and $\Omega=\times_{i\in \mathcal{N}}\Omega_i$.
Write algorithm \eqref{gcd} in a compact form  as
\begin{equation}
\begin{array}{lll}
\dot{X}       & = & P_{\Omega}( X-  \nabla f(X) +\Lambda)-X,\\
\dot{\Lambda} & = & -(L\otimes I_m) \Lambda-(L\otimes I_m) Z+(D-X),\\
\dot{Z}       & = & (L\otimes I_m) \Lambda.  \label{cgcd}
\end{array}
\end{equation}

\begin{thm}\label{correct}
Under Assumptions \ref{asumFun}-\ref{assum2},
if the initial point $x_i(0) \in \Omega_i, \;  \forall i \in \mathcal{N}$,  then $x_i(t) \in \Omega_i,\; \forall t \geq 0,\forall i \in \mathcal{N}$, and
$col(X^*,\Lambda^*,Z^*)$ is the equilibrium point of the distributed algorithm \eqref{gcd} with $X^{*}$ as  the optimal solution to \eqref{CRA1}.
\end{thm}

{\bf Proof}:  Note that
$\dot{x}_i \in T_{\Omega_i}(x_i), \; \forall x_i\in \Omega_i$ because $P_{\Omega_i}(x_i-\nabla f_i(x_i)+\lambda_i) \in \Omega_i$.
Given the initial point $x_i(0) \in \Omega_i$, by  Nagumo's theorem (referring to page 174 and page 214 of \cite{aubin}),  $x_i(t)\in \Omega_i, \forall t\geq 0$ (Related proofs can also be found in \cite{arrow2}).

To obtain the equilibrium point, we set $\dot{Z}=\mathbf{0}_{mn}$ and get $\Lambda^*= \mathbf{1}_n \otimes  \lambda^*$, $\lambda^*\in \mathbf{R}^m $,
because the graph $\mathcal{G}$ is connected.
$\dot{\Lambda}=\mathbf{0}_{mn}$ implies that  $ (L\otimes I_m) Z^* = D-X^*$. Because the graph $\mathcal{G}$ is undirected,
 $\mathbf{1}_n^T L=\mathbf{0}_n^T$ and  $(\mathbf{1}^T_n \otimes I_m) (L\otimes I_m) Z = (\mathbf{1}^T_n  L)\otimes I_m Z =(\mathbf{1}_n^T \otimes I_m)(D-X)=\mathbf{0}_{m}$.
Hence,  $\sum_{i\in \mathcal{N}} d_i=\sum_{i\in \mathcal{N}} x_i^*$.
Then, at the equilibrium point,
$$\Lambda^* =\mathbf{1}_n\otimes \lambda^*, \; \lambda^* \in \mathbf{R}^m, (L \otimes I_m) Z^*= D-X^*, \;  \sum_{i\in \mathcal{N}} d_i =\sum_{i\in \mathcal{N}}x_i^*.  $$
Also, at the equilibrium point, $\dot{x_i}=0$ implies that
$P_{\Omega_i}( x^*_i-  \nabla f_i(x^*_i) + \lambda^*)=x^*_i$.
It follows from Lemma 2.38 of \cite{ob} that
$$ - \nabla f_i(x^*_i) + \lambda^* \in C_{\Omega_i}(x^*_i), i=1,...,n . $$

Therefore, the equilibrium point  $col(X^*,\Lambda^*,Z^*)$ of  \eqref{gcd} satisfies
\begin{equation}\label{kkt2}
\begin{array}{lll}
&&\mathbf{0}_{mn} \in   \nabla f(X^*)  - \mathbf{1}_n \otimes \lambda^*  +  C_{\Omega}(X^*), \\
&& (\mathbf{1}^T_{n} \otimes I_m)  X^*= (\mathbf{1}^T_{n}\otimes I_m) D, \; \quad X^* \in \Omega,
 \end{array}
\end{equation}
which is exactly the optimality condition (KKT) for DRAO \eqref{CRA1} by Theorem 3.34 in \cite{ob}.
Thus, the conclusion follows.
\hfill $\Box$

\begin{rem}
It can be shown that the equilibrium point of \eqref{gcd} has
$\lambda^*$ as the dual optimal solution to  problem \eqref{problem2},
and $Z^*$ as the dual optimal solution to problem \eqref{problem3},
following a similar analysis routine of Theorem \ref{correct} and Proposition 5.3.2 in \cite{bersekas}.
It can also be shown that any $col(X^*, \mathbf{1}_m\otimes \lambda^*, Z^*)$, with $X^*$ as the optimal solution to \eqref{CRA1}, $\lambda^*$ as the
dual optimal solution to \eqref{problem2} and $Z^*$ as the dual optimal solution to \eqref{problem3}, corresponds to an equilibrium point of \eqref{gcd}.
We do not discuss their details here for space limitations.
\end{rem}

\begin{rem}\label{kktcompare}
The differences between our work and some existing ones are listed as follows:
\begin{itemize}
\item The KKT condition for DRAO without LFC is
\begin{equation}\label{kkt3}
\begin{array}{lll}
\nabla f(X^*)  = \mathbf{1}_n \otimes \lambda^*,  \;  (\mathbf{1}^T_{n} \otimes I_m)  X^*= (\mathbf{1}^T_{n}\otimes I_m) D.
\end{array}
\end{equation}
The KKT condition  \eqref{kkt2} for DRAO with LFC and the condition  \eqref{kkt3} for DRAO without LFC  are  totally different.
\eqref{kkt3} requires the optimal allocations to be the  points  satisfying the network resource constraint with the same marginal costs (gradients), while
\eqref{kkt2} requires the  optimal allocations  to be  feasible in both network resource constraint and LFCs.
The optimal allocations in \eqref{kkt2} should also satisfy a variational inequality related to both the objective functions' gradients and the normal cones of the LFC sets.
In fact, the marginal costs (gradients) at the optimal allocations of \eqref{kkt2}  do  not necessarily reach the same levels,  and  the differences can be seen as the ``price of allocation feasibility".

\item
The previous algorithms (except the one in \cite{EDcc}) kept the network resource constraint satisfied (ensured its eventual feasibility)  by setting feasible initial points through  the initialization coordination procedure.
In other words, the network resource constraint can be guaranteed {\bf only if} it is satisfied at the initial moment.
However,  the initialization for the network resource constraint  is quite restrictive for large-scale dynamical networks because it involves global coordination
 and has to be performed  every time the network data/configuration changes.
 Moreover,  it is not trivial to achieve  the initialization coordination with both the  LFCs and the network resource constraint
 (refer to \cite{ED3} for an initialization procedure with one dimensional interval constraint).

\item  In fact,  proportional-integral (PI)-type consensus dynamics (\cite{wang} and  \cite{cort2}) and projected gradient flows (\cite{liu}), which both have been utilized for distributed optimization,
are  combined together to obtain the algorithm \eqref{gcd} for the KKT condition \eqref{kkt2}.
  Local $\lambda_i$ acts as the local shadow price,
  and all the local shadow prices  must reach consensus to be the global market clearing price.
  Therefore, a second-order PI-based consensus dynamics is incorporated into (\ref{gcd})
  with the integral variable $z_i$ summing up the disagreements between $\lambda_i$ and $\{\lambda_j, j \in \mathcal{N}_i\}$.
  Meanwhile, the dynamics of $x_i$  adjusts the local  allocation by comparing the local shadow price and the local gradient,
  and also utilizes the  projection operation in order to make the local allocation flowing  within its LFC set all the time.
\end{itemize}
\end{rem}

Therefore, the algorithms given in \cite{RA3},  \cite{RA5}, \cite{RA6} and \cite{RA7} failed to solve \eqref{CRA1} because
they cannot ensure the LFCs given in \eqref{kkt2}.
Note that the algorithm \eqref{gcd} can  ensure LFCs even  during the algorithm flow with projection operations.
It only requires that  each agent has its initial  allocation belonging to  its LFC set,
which can be trivially accomplished by each agent with one-step  local projection operation.
Furthermore, algorithm \eqref{gcd}  ensures the network  resource constraint
asymptotically without caring about whether it is satisfied at the initial points, and therefore,
is free of initialization coordination procedure.
Due to free of any center and  initialization,
algorithm \eqref{gcd} can adaptively handle online  data  without  re-initialization whenever the network data/configuration changes.
Moreover,  it can work in a ``plug-and-play" manner  for dynamical networks with leaving-off or plugging-in of agents.

Next, let us analyze the convergence of \eqref{gcd},  The analysis techniques are inspired by the  projected dynamical systems for variational inequalities (referring to \cite{jun} and \cite{Gao}) and distributed optimization (referring to \cite{shi1}, \cite{cor1}, \cite{xie} and  \cite{liu}).

\begin{thm}\label{ascgcd}
Under Assumptions \ref{asumFun}-\ref{assum2}, and given bounded initial points $x_i(0) \in \Omega_i, \;  \forall i \in \mathcal{N}$, the trajectories of the algorithm \eqref{gcd} are bounded and  converge to an  equilibrium point of \eqref{gcd}, namely,
agent $i$ asymptotically achieves its optimal allocation $x_i^*$ of \eqref{CRA1} with \eqref{gcd}.
\end{thm}

{\bf Proof}:
Take $m=1$ without loss of generality.
Denote  $\bar{\Omega}=\Omega \times \mathbf{R}^{n}\times \mathbf{R}^{n} .$
Define a new vector $S=col(X,\Lambda,Z)$  and the vector function $F(S): R^{3n}\rightarrow R^{3n}$ as
\begin{equation}
F(S)=\left(
       \begin{array}{c}
         \nabla f(X)-\Lambda \\
         L\Lambda+LZ-(D-X) \\
         -L\Lambda \\
       \end{array}
     \right). \nonumber
\end{equation}

Recalling the form in \eqref{cgcd} and the fact that
$P_{\mathbf{R}^{n}}(x)=x,\; \forall x\in \mathbf{R}^n$, the dynamics of
all the agents can be written as
$$
\dot{S} = P_{\bar{\Omega}}(S-F(S))-S.
$$
Define $H(S)=P_{\bar{\Omega}}(S-F(S))$, and then give a Lyapunov function as
\begin{equation}
V_{g} = -\langle F(S), H(S)-S\rangle -\frac{1}{2} || H(S)-S ||_2^2+ \frac{1}{2}  || S-S^* ||_2^2, \nonumber
\end{equation}
where $S^*=col(X^*,\Lambda^*,Z^*)$, $X^*$ is the optimal solution to  \eqref{CRA1},
$\Lambda^*$ is the optimal solution to \eqref{problem3}, and $Z^*$ is the dual optimal solution to \eqref{problem4}.

Notice that $X^*$ is  a finite point  from Assumption \ref{slater}.
Because the Slater's condition is satisfied with Assumption \ref{slater},
the dual optimal solution $\lambda^*$  for \eqref{problem2}  exists and is finite by
Proposition 5.3.1 of \cite{bersekas}.
Since the function $f_i(x_i)$ is strictly convex, $\nabla q_i(\lambda^*) = d_i-x_i^*$ by Theorem 2.87 of \cite{ob} and the saddle point property of \eqref{problem2}.
Then the KKT condition of \eqref{problem3} is $L\Lambda^*=\mathbf{0}$ and  $D-X^*-L\Lambda^*-LZ^*=\mathbf{0}$. Hence, $LZ^*=D-X^*$ implies the finiteness of the dual optimal solution $Z^*$ for \eqref{problem3}.
Therefore, $S^*$ is a finite point.

In fact, with the KKT conditions to \eqref{problem2} and \eqref{problem3},
$ F(S^*)=col(   (\nabla f(X^*) - \Lambda^*), \mathbf{0}_n, \mathbf{0}_n), $
  and $-\nabla f(X^*) + \Lambda^*\in C_{\Omega}(X^*)$, we have
\begin{equation}\label{equaligcd}
 H(S^*) = P_{\bar{\Omega}}(S^*-F(S^*))=S^*, \;  - F(S^*) \in  C_{\bar{\Omega}}(S^*).
\end{equation}
Due to \eqref{projection} and \eqref{projection2},
\begin{equation}
\begin{array}{lll}
&&-\langle F(S), H(S)-S\rangle - \frac{1}{2}\langle H(S)-S, H(S)-S \rangle \\
&&= \frac{1}{2} [ ||F(S) ||_2^2- ||F(s)+H(S)-S ||_2^2]\\
&&= \frac{1}{2}[ ||S- F(S)-S ||_2^2 - ||H(S)-(S-F(S)) ||_2^2 ]\\
&&\geq \frac{1}{2} ||S-H(S) ||_2^2.  \nonumber
\end{array}
\end{equation}

Hence, $V_g\geq \frac{1}{2} ||S-H(S)||_2^2 +\frac{1}{2} ||S-S^* || \geq 0, $ and
$V_g=0$ if and only if $S=S^*$.

By Theorem 3.2 of \cite{fukusma}, any asymmetric variational inequality can be converted to a differentiable optimization problem. As a result,
\begin{equation}
\dot{V}_g = (F(S)-[J_F(S)-I](H(S)-S)+ S-S^*)^T (H(S)-S), \nonumber
\end{equation}
where $J_F(S)$ is the Jacabian matrix of $F(S)$ defined as
\begin{equation}
J_F(S)=  \left(
  \begin{array}{ccc}
    \nabla^2 f(X) & -I & 0 \\
    I & L & L \\
    0 & -L & 0 \\
  \end{array}
\right). \nonumber
\end{equation}

With Assumptions \ref{asumFun} and \ref{assum2},
$$
S^T J_F(\bar{S}) S = X^T \nabla^2 f(\bar{X}) X + \Lambda^T L \Lambda > 0,  \quad \forall \bar{S} \in \bar{\Omega}, \; \forall S \neq \mathbf{0}\in \mathbf{R}^{3n}.$$

With \eqref{projection}, taking $x=S-F(S)$ and $y=S^*$ gives
$\langle S-F(S)-H(S), H(S)-S^* \rangle \geq 0$, which implies
 $\langle S-H(S)-F(S), H(S)-S+S-S^* \rangle \geq 0.$
Hence,
 $- ||H(S)-S ||_2^2 + \langle S-H(S),   S-S^* \rangle  + \langle -F(S),  H(S)-S\rangle + \langle -F(S), S-S^*   \rangle \geq
 0 $, or equivalently,
 \begin{equation}\label{proest}
  \langle S-H(S),    S-S^*+F(S)  \rangle  \geq ||H(S) -S||_2^2  +  \langle F(S), S-S^*   \rangle. \nonumber
  \end{equation}
Consequently,
\begin{equation}
\begin{array}{lll}
\dot{V_g} &\leq  & - (H(S)-S)^T J_F(S) (H(S)-S) + || H(S)-S||_2^2 \\
         & \; &+    \langle S-S^*+F(S), H(S)-S  \rangle  \\
         & \leq  & -\langle F(S), S-S^*   \rangle \\
         & \leq & - \langle F(S), S-S^*    \rangle + \langle F(S^*), S-S^*    \rangle - \langle F(S^*), S-S^*    \rangle \\
           & \leq & - \langle  F(S)-F(S^*), S-S^*  \rangle  - \langle F(S^*), S-S^*    \rangle. \nonumber
\end{array}
\end{equation}
In fact,
$\langle  F(S)-F(S^{'}), S-S^{'}\rangle
 =  \langle \nabla f(X)-\nabla f(X^{'}),X-X^{'} \rangle + \langle \Lambda-\Lambda^{'}, L(\Lambda-\Lambda^{'})\rangle
 \geq 0, \; \forall S,\;S^{'}\in \bar{\Omega},$
because the local objective functions are convex, and the Laplacian matrix is positive semi-definite by Assumptions \ref{asumFun} and \ref{assum2}.
With \eqref{equaligcd}, we have  $ \langle F(S^*), S-S^*   \rangle \geq 0 $.
Then $\dot{V_g} \leq 0,$ and any finite equilibrium point $S^*$ of \eqref{gcd} is Lyapunov stable.
Furthermore, there exists a  forward compact  invariance set given any finite initial points,
\begin{equation}\label{compactset}
I_S=\{col(X,\Lambda,Z)\big{|}\frac{1}{2}{||S-S^* ||}\leq V_g(S(0)) \}.
\end{equation}

Therefore, with the local Lipstchitz continuity of  the objective functions' gradients in Assumption \ref{asumFun} and
the non-expansive property of projection operation \eqref{projection3},
$P(S-F(S))-S$ is Lipstchitz over the compact set $I_S$ in \eqref{compactset}.
There exists a unique solution to \eqref{gcd} with time domain $[0,\infty).$
Also, the compactness and convexity of $I_S$
implies the existence of equilibrium point to dynamics \eqref{gcd} (referring to page 228 of \cite{aubin}).
By Theorem \ref{correct} and without loss of generality, the equilibrium point is assumed to be  $S^*$.

Furthermore, there exists $c^*\in  C_{\Omega}(X^*)$ such that  $ -\nabla f(X^*)+\Lambda^* =c^*$, $LZ^*=D-X^*$, and $\Lambda^*= \mathbf{1}_n\lambda^*$.
\begin{equation}
\begin{array}{lll}
\dot{V_g} &\leq &   -  \langle F(S), S-S^*   \rangle = -\langle X-X^*, \nabla f(X)-\Lambda\rangle\\
                 && -  \langle \Lambda-\Lambda^*, L\Lambda+LZ-(D-X^*)\rangle- \langle Z-Z^*, -L\Lambda\rangle\\
          & \leq &  -  \langle X-X^*, \nabla f(X)-\Lambda -\nabla f(X^*)+\Lambda^* -c^*  \rangle  \\
                 && -  \langle \Lambda-\Lambda^*, L(\Lambda-\Lambda^*) \rangle
                     -  \langle \Lambda-\Lambda^*, L(Z-Z^*) \rangle\\
                 &&  -  \langle \Lambda-\Lambda^*, LZ^*-(D-X)\rangle
                  -  \langle Z-Z^*, -L(\Lambda-\Lambda^*)\rangle \\
          &\leq &   - \langle X-X^*, \nabla f(X)-\nabla f(X^*)  \rangle  +  \langle X-X^*,  c^*  \rangle \\
                 && - \langle \Lambda-\Lambda^*, L(\Lambda-\Lambda^*) \rangle\\
          & \leq &  - \langle X-X^*, \nabla f(X)-\nabla f(X^*)  \rangle
           -\langle \Lambda-\Lambda^*, L(\Lambda-\Lambda^*) \rangle,  \nonumber
\end{array}
\end{equation}
where the last step follows from $ \langle X-X^*,  c^*  \rangle \leq  0$.
Denote the set of points satisfying $\dot{V}_g=0$ as $E_g=\{(X,\Lambda,Z) \big |  \dot{V}_g=0 \}.$
Because the Hessian matrix of $\nabla^2 f(X)$ is positive definite,
$
\nabla f(X) = \nabla f(X^*)+ \int_{0}^1 \nabla^2 f(\tau X+(1-\tau)X^*)^T(X-X^*)d\tau
$
and the null space for $L$ imply that
$$E_g=\{(X,\Lambda,Z) \big | X=X^*,\Lambda \in span\{\alpha \mathbf{1}_n\}  \}.$$

Then we claim that  the maximal invariance set within the set $E_g$ is exactly the equilibrium point of \eqref{gcd}.
In fact, $\Lambda\in span\{\alpha \mathbf{1}_n\} $ implies $Z=Z^*$.  Hence, $\dot{\Lambda}=LZ^*-(D-X^*)$. However,
$LZ^*-(D-X^*)$ must be zero; otherwise $\Lambda$ will go to infinity, which contradicts that $E_g$ is a compact set within $I_S$.
Thus, $\dot{\Lambda}=0$ and $\Lambda=\Lambda^*$.

By the LaSalle invariance principle and Lyapunov stability of the equilibrium point,  the system \eqref{gcd} converges to its equilibrium point, which implies the conclusion.
\hfill $\Box$

\begin{rem}
Although an initialization-free algorithm has also been proposed and
investigated for the DEDP, a special case of DRAO,  in \cite{EDcc},
algorithm \eqref{gcd} provides a  different  algorithm to
address the DRAO problem without initialization.
Additionally, algorithm \eqref{gcd} can handle general
multi-dimensional LFCs explicitly with the projection operation, while
\cite{EDcc} only addressed one-dimensional box constraints
 with a penalty method.  Moreover, one agent was
required to know the total network resource all the time in a
time-varying resource case in \cite{EDcc},  while each agent only
knows its local resource  in \eqref{gcd}. Moreover, our
techniques introduce a variational inequality viewpoint in
addition to Lyapunov methods with the invariance principle.
\end{rem}

The following  example illustrates  how \eqref{gcd} ``adaptively" achieves the  optimal resource allocation
without re-initialization for a dynamical network.
Notice that the following example cannot be directly addressed by the algorithm in \cite{EDcc}.

\begin{table}\label{tp1}
\caption{Parameters setting of Example \ref{exa1}}
\label{Simulation result}
\begin{center}
\begin{tabular}{  | l | l | l | l  | }
\hline
\quad         & $0\sim600s $       & $ 600s\sim1200s   $             & $1200s\sim $            \\ \hline
$a_1,d_1$     &  (8,2),    (8,2)   &  (0.1,0.3),    (8,2)           & (0.1,0.3),    (12,-3)               \\ \hline
$a_2,d_2$     &  (4,7),    (3,4)   &  (-17,15),     (3,4)           & (-17,15),     (0,7)            \\ \hline
$a_3,d_3$     &  (0.13,8), (3,8)   &  (0.13,8),     (-5,12)         & (3,0.7),      (-5,12)           \\ \hline
$a_4,d_4$     &  (4,20),   (10,2)  &  (4,20),       (1,15)          & (5,17),       (1,15)             \\ \hline
\end{tabular}
\end{center}
\end{table}

\begin{figure}
  \centering
  \includegraphics[width=3in]{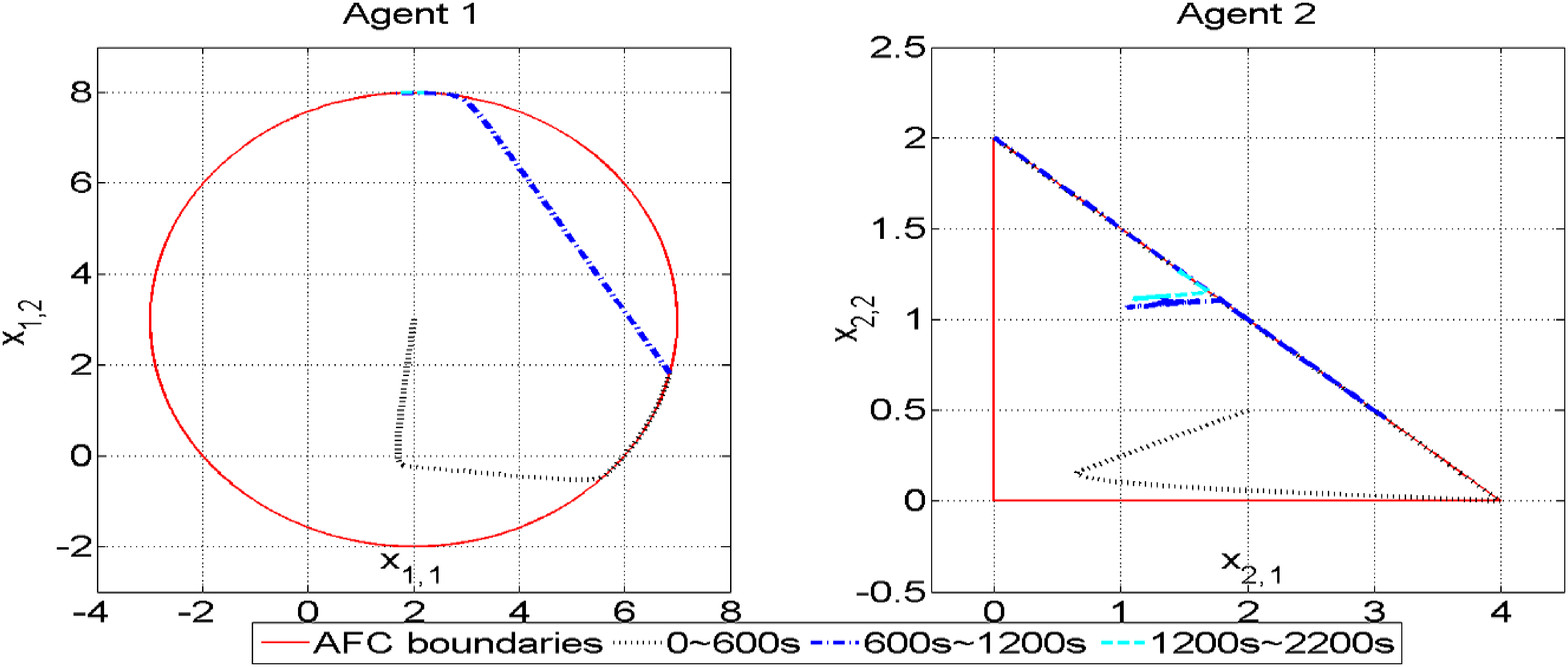}\\
  \caption{\small{The trajectories of the allocations of agent $1$ and agent $2$.}}\label{gcds1}
\end{figure}

\begin{figure}
  \centering
  \includegraphics[width=3in]{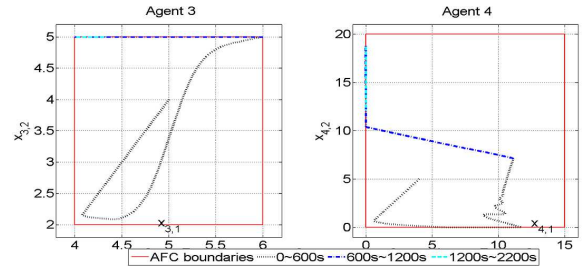}\\
  \caption{\small{The trajectories of the allocations of agent $3$ and agent $4$.}}\label{gcds2}
\end{figure}

\begin{exa}\label{exa1}
Four agents cooperatively optimize  problem \eqref{CRA1}.
The allocation variable and resource data for agent $i$ are  $x_i=(x_{i,1},x_{i,2})^T \in \mathbf{R}^2,$ $d_i=(d_{i,1},d_{i,2})^T \in \mathbf{R}^2$, respectively.
The objective functions $f_i(x_i)$ are  parameterized with  $a_i=(a_{i,1},a_{i,2})^T \in \mathbf{R}^2$  as follows:
\begin{equation}
f_i(x_i)  =  (x_{i,1}+a_{i,1}x_{i,2})^2 + x_{i,1}+ a_{i,2}x_{i,2}+0.001(x^2_{i,1}+x^2_{i,2}). \nonumber
\end{equation}

The LFCs of the four agents are given as
follows: $\Omega_1=\{x_1\in \mathbf{R}^2| (x_{1,1}-2)^2+ (x_{1,2}-3)^2 \leq
25\}$, $\Omega_2=\{x_{2}\in \mathbf{R}^2 | x_{2,1}\geq 0,
x_{2,1}\geq 0, x_{2,1}+2x_{2,2}\leq 4\}$, $ \Omega_3=\{x_3 \in
\mathbf{R}^2 | 4 \leq x_{3,1}\leq 6, 2\leq x_{3,2}\leq 5\}$ and
$\Omega_4=\{x_4\in \mathbf{R}^2 |  0\leq x_{4,1}\leq 15, 0\leq
x_{4,2}\leq 20\}$, respectively, and their boundaries are  shown in
Figures \ref{gcds1} and \ref{gcds2}.

The agents share information with a ring graph $\mathcal{G}$:
$$1\leftrightarrow 2 \leftrightarrow 3 \leftrightarrow  4 \leftrightarrow 1.$$
The initial allocation  $x_i(0)$ of agent $i$ in algorithm
(\ref{gcd}) is randomly chosen within its  LFC  set, and $\Lambda,
Z$ are set with zero initial values. The data $a_i$ and $d_i$
switches as   Table \ref{tp1}, while the allocation variables remain
unchanged  when  the data switches. The simulation results are shown
in Figures \ref{gcds1}, \ref{gcds2}, \ref{gcds3}, and \ref{gcds4}.

\begin{figure}
  \centering
  \includegraphics[width=3in]{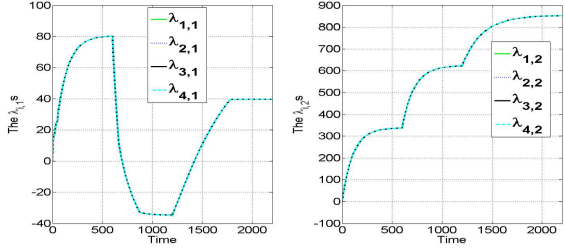}\\
  \caption{ \small{The trajectories of the  Lagrangian multiplies $\lambda_{i}$'s }}\label{gcds3}
\end{figure}

\begin{figure}
  \centering
  \includegraphics[width=3in]{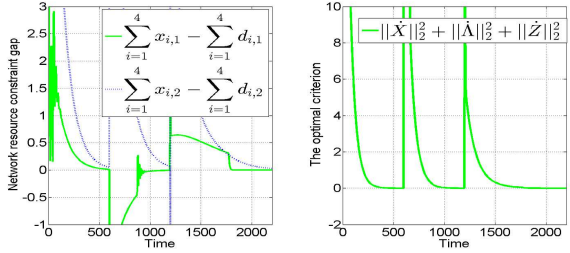}\\
  \caption{ \small{Network resource constraint and optimality condition}}\label{gcds4}
\end{figure}

Figures \ref{gcds1} and \ref{gcds2} show that the
agents' allocation variables  always remain within the corresponding
LFC sets, while Figure \ref{gcds3} shows that the Lagrangian
multipliers reach consensus after the transient processes. Figure
\ref{gcds4} demonstrates that the network resource constraint  can
be satisfied asymptotically even though it is violated each time the
data/configuration changes, and Figure \ref{gcds4} also reveals that
$|| \dot{X}||_2^2+||\dot{\Lambda} ||_2^2+ || \dot{Z}||_2^2$ always
converges to zero,  guaranteeing the optimality of the resource
allocations.
\end{exa}

\section{Differentiated projected algorithm for DRAO}

In this section,  the differentiated  projection operator \eqref{dpoperator} is applied to
derive an  algorithm for \eqref{CRA1}. 
In fact, the projected dynamics based on the operator in
$\eqref{dpoperator}$ was firstly introduced in the study of
constrained stochastic approximation in \cite{proddup2}, and later
was utilized to solve variational inequalities and constrained
optimization problems in \cite{prods2}, \cite{prods4}, and
\cite{cortes_prods}. Here the operator \eqref{dpoperator} is
applied to the construction  of the distributed resource allocation
algorithm for agent $i$ given as follows:
\begin{equation}\label{cd}
\begin{array}{l}\hline
{\bf Differentiated \; projected \;  algorithm \; for \; agent \; i}:           \\ \hline
\displaystyle  \dot{x}_i       = \Pi_{\Omega_i}(x_i, -\nabla f_i(x_i)+\lambda_i)            \\
\displaystyle  \dot{\lambda}_i = -\sum_{j\in \mathcal{N}_i}(\lambda_i-\lambda_j)-\sum_{j\in \mathcal{N}_i}(z_i-z_j)+(d_i-x_i)  \\
\displaystyle  \dot{z_i}       = \sum_{j\in \mathcal{N}_i}(\lambda_i-\lambda_j)                              \\  \hline
\end{array}
\end{equation}

\begin{rem}
The algorithm \eqref{cd} is a direct extension of \eqref{gcd} by differentiating the projection operator, where each agent is required to project $-\nabla f_i(x_i)+\lambda_i$ onto the tangent cone $T_{\Omega_i}(x_i)$.
Thereby, \eqref{cd} has the additional burden for the tangent cone computation compared with \eqref{gcd}.
However, for some specific  convex sets such as polyhedron, Euclidean ball, and boxes, it is not hard to
get the close form of the tangent cone  at any given point.
\end{rem}

Similar to the algorithm \eqref{gcd}, the algorithm \eqref{cd} is also a distributed algorithm, and
it does not need any initialization coordination procedure.  Therefore, it
can efficiently process online data for  dynamical networks.

Although algorithm \eqref{cd} is a discontinuous
dynamical system, the solution to \eqref{cd} is
well-defined in the Caratheodory sense (an absolutely continuous
function $col(X(t),\Lambda(t),Z(t)): [0,T]\rightarrow \mathbf{R}^{3mn}$
is a solution of \eqref{cd} if \eqref{cd}  is satisfied for almost
all $t\in [0,T]$, referring to Definition 2.5 in \cite{prods2}). The
existence of an absolutely continuous solution to \eqref{cd}  can be
found in Theorem 3.1 of \cite{prods3}, and the condition when the
solution can  be extended to interval $[0,\infty]$ is given in
Theorem 1 of \cite{prods4}.

The following result shows the correctness of algorithm \eqref{cd}.

\begin{thm}
Suppose that Assumptions \ref{asumFun}-\ref{assum2} hold.
If the initial point $x_i(0) \in \Omega_i, \; \forall i \in \mathcal{N}$,  then $x_i(t) \in \Omega_i, \; \forall t \geq 0,\forall i \in \mathcal{N}$,  and there is the equilibrium point of the algorithm \eqref{cd} with $X^{*}=col(x_1^*,...,x_n^*)$  as  the optimal solution to \eqref{CRA1}.
\end{thm}

{\bf Proof}: Obviously,  $\dot{x}_i \in
T_{\Omega_i}(x_i), \; \forall x_i\in \Omega_i$ according to Lemma
\ref{dpoperatorbaic}. It follows that \eqref{cd} has an absolutely
continuous solution on interval $[0,\infty]$ by Theorem 1 of
\cite{prods4}. Moreover,  Theorem 3.2 of \cite{prods3} shows that
the solution of  \eqref{cd} coincides with a slow solution of a
differential inclusion.  Given the
initial point $x_i(0)\in \Omega_i$, $x_i(t)\in \Omega_i, \forall
t\geq 0$ holds in light of the viability theorem in \cite{aubin}.

By Lemma \ref{dpoperatorbaic}, we have that $\Pi_{\Omega_i}(x_i,-\nabla f_i(x_i)+\lambda_i)=0$ if
at least one of the following cases is satisfied:
(i):$x_i\in int(\Omega_i),$ and $-\nabla f_i(x_i)+\lambda_i=0$;
(ii): $x_i\in \partial \Omega_i$, and $-\nabla f_i(x_i)+\lambda_i=0$;
(iii): $x_i\in \partial \Omega_i$, and $-\nabla f_i(x_i)+\lambda_i\in C_{\Omega_i}(x_i)$.
Hence,  $\Pi_{\Omega_i}(x_i,-\nabla f_i(x_i)+\lambda_i)=0$ implies $-\nabla f_i(x_i)+\lambda_i\in C_{\Omega_i}(x_i)$.
Following similar analysis of Theorem \ref{correct}, at the equilibrium point we have
\begin{equation}
\begin{array}{lll}\label{eq2}
&& \Lambda^*=\mathbf{1}_n \otimes \lambda^*,\lambda^*\in \mathbf{R}^m, \quad \; x^*_i \in \Omega_i \\
&& (L\otimes I_m) Z^* = D-X^*, \sum_{i\in \mathcal{N}} x^*_i  =  \sum_{i\in \mathcal{N}} d_i, \\
&& -\nabla f_i(x_i^*)+ \lambda^*  \in C_{\Omega_i}(x_i^*).
\end{array}
\end{equation}
Thus, the optimality condition \eqref{kkt2} of \eqref{CRA1} is satisfied by the equilibrium point of \eqref{cd}.
\hfill $\Box$

Next result shows the convergence of \eqref{cd} when the local objective functions are strongly convex.

\begin{thm}
Suppose that Assumptions \ref{asumFun}-\ref{assum2} hold, and
the local objective functions $f_i(x_i)$ are $\mu_i$-strongly convex functions with
$k_i$-Lipschitz continuous gradients.
Given bounded  initial points $x_i(0) \in \Omega_i, \; \forall i \in \mathcal{N}$,
the trajectories of  algorithm (\ref{cd}) converge to its equilibrium point.
Furthermore, if there are no LFCs (that is, $\Omega_i=\mathbf{R}^m,i=1,\cdots,n$), then algorithm (\ref{cd})
exponentially converges to its equilibrium point.
\end{thm}

{\bf Proof}:  We still take $m=1$ without loss of generality.

At first, we show the convergence of \eqref{cd}.
By Lemma \ref{dpoperatorbaic},
\begin{equation}\label{eq1}
\Pi_{\Omega_i} (x_i, -\nabla f_i(x_i) + \lambda_i) =  -\nabla f_i(x_i) + \lambda_i - \beta(x_i)n_{i}(x_i),
\end{equation}
where $ n_i(x_i)\in c_{\Omega_i}(x_i),\; \beta(x_i)\geq 0$.
Notice that there exist $\beta(x^*_i)\geq 0$ and $n_{i}(x^*_i)\in c_{\Omega_i}(x^*_i)$ at the equilibrium point such that
$
\nabla f_i(x^*_i) - \lambda^* =  -\beta(x^*_i)n_{i}(x^*_i).
$

Define the following variables
\begin{equation}
\begin{array}{lll}\label{corc}
Y=X-X^*,  \;              & \theta=[r \; R]^T  Y,  &  Y= [r \; R]\theta,\\
V=\Lambda-\Lambda^*, \;    & \eta=[r \; R]^T    V , &  V=[r \; R] \eta,\\
W=Z-Z^*,                  & \delta= [r \;  R]^T W,  &  W=[r \; R]\delta,\\
\end{array}
\end{equation}
with
 $r=\frac{1}{\sqrt{n}} \mathbf{1}_n$ and $r^T R=\mathbf{0}^T_n$ such that  $R^T R=I_{n-1}$ and $R R^T = I_n- \frac{1}{n} \mathbf{1}_n\mathbf{1}_n^T$.
We partition the variables $\theta,\eta,\delta$ as $col(\theta_1,\theta_2)$, $col(\eta_1,\eta_2)$, $col(\delta_1,\delta_2)$ with $\theta_1,\eta_1,\delta_1 \in \mathbf{R}$ and
$\theta_2,\eta_2,\delta_2 \in \mathbf{R}^{n-1}$.

Then the dynamics of the  variables $\theta,\eta,\delta$
can be derived with \eqref{cd}, \eqref{eq2}, \eqref{eq1} and \eqref{corc} as follows,
\begin{equation}\label{corcdy}
\begin{array}{lll}
\dot{\theta}_1  =    -r^T h+\eta_1; & \; &  \dot{\theta}_2  =    -R^T h +\eta_2;\\
\dot{\eta}_1    =    -\theta_1    ; & \; &  \dot{\eta}_2    =    -\theta_2 -R^T L R\eta_2 - R^T L R \delta_2;\\
\dot{\delta}_1  =    0            ; & \; &  \dot{\delta}_2  =    R^TLR \eta_2,
\end{array}
\end{equation}
where $h=\nabla f(Y+X^*)-\nabla f(X^*) + N_{\Omega}(X)-N_{\Omega}(X^*)$,
$N_{\Omega}(X) = col(\beta(x_1)n_{1}(x_1),..., \beta(x_n)n_n(x_n)     ), $ and
$N_{\Omega}({X}^*)= col( \beta({x}^*_1)n_{1}({x}^*_1),...., \beta({x}^*_n)n_{n}({x}^*_n)   )$.

Construct the following function
\begin{equation}\label{vs1}
\begin{array}{lll}
 V^s_1 & =  & \frac{1}{2}\alpha (\theta^T \theta+ \eta^T\eta) + \frac{1}{2}(\alpha+\gamma) \delta^T_2 \delta_2  \\
   & +  & \frac{1}{2} \gamma  (\eta_2+ \delta_2 )^T (\eta_2+\delta_2 ),
\end{array}
\end{equation}
where $\alpha,\gamma$ are positive constants to be determined later.
Obviously,
$\frac{1}{2} \alpha ||p ||_2^2 \leq V^s_1 \leq (\frac{1}{2}(\alpha+\gamma)+\gamma) ||p ||_2^2$ where $p=col(\theta,\eta,\delta_2)$.
The derivative of $V_1^s$ along (\ref{cd}) is
\begin{equation}
\begin{array}{lll}
\dot{V}^s_1 & = &  \alpha(- Y^T h -\eta_2 R^T L R \eta_2 -\eta_2 R^T L R \delta_2)\\
        & + &  \gamma (-\eta_2^T  \theta_2- \eta_2 R^T L R \delta_2 - \delta_2 \theta_2 - \delta_2^T R^T L R \delta_2)\\
        & + &  (\alpha+\gamma)\delta_2^T R^T L R \eta_2. \nonumber
\end{array}
\end{equation}

Because $n_i(x_i^*)\in c_{\Omega_i}(x_i^*)$ and $ \beta(x_i^*)\geq 0$, $ \beta{(x_i^*)} \langle x_i-x_i^*,  n_i(x_i^*)   \rangle  \leq  0$.
Moreover, $n_i(x_i)\in c_{\Omega_i}(x_i)$ and $\beta(x_i)\geq 0$ imply  that  $\beta(x_i) \langle x_i-x_i^*,n_i(x_i)\rangle \geq 0 $.
Because the local objective functions are strongly convex,
$
 Y^T h =   (X-X^*)^T(  \nabla f(X)-\nabla f(X^*) + N_{\Omega}( X)- N_{\Omega}(X^*))
        =   (X-X^*)^T(  \nabla f(X)-\nabla f(X^*) )
       + \sum_{i=1}^n \langle x_i-x_i^*,+ \beta{(x_i)}n_{i}(x_i) \rangle
        +  \sum_{i=1}^n  \langle x_i-x_i^*,  -\beta(x^*_i)n_i(x^*_i)\rangle
       \geq \bar{\mu} \theta^T \theta
$
 where $\bar{\mu}=\min\{\mu_1,\cdots,\mu_n\}.$

Denote $s_1\leq s_2\leq ...\leq s_n$ as the ordered eigenvalues of Laplacian matrix $L$.
Obviously, $s_1=0$ and $s_2>0$ when the graph $\mathcal{G}$ is connected.
 By \eqref{lap},
$$\dot{V}^s_1   \leq      -\alpha\bar{\mu}\theta^T \theta -\alpha s_2 \eta_2^T \eta_2  -\gamma s_2 \delta_2^T \delta_2
              - \gamma\eta_2^T  \theta_2 - \gamma\delta_2  \theta_2.$$

According to
$-\gamma \eta_2^T \theta_2\leq \frac{1}{2}\gamma^2 \theta_2^T \theta_2 + \frac{1}{2} \eta_2^T \eta_2,$
and
$ - \gamma\delta_2^T \theta_2 \leq \frac{1}{2} \gamma^2 \theta_2^T \theta_2 + \frac{1}{2}\delta_2^T \delta_2,$
we have
\begin{equation}\label{v1dot}
\dot{V}^s_1   \leq     -(\alpha \bar{\mu} -\gamma^2) \theta^T \theta - (\alpha s_2 - \frac{1}{2}) \eta_2^T\eta_2
           -      (\gamma s_2 -\frac{1}{2}) \delta_2^T \delta_2.
\end{equation}
Take $\gamma$ and $\alpha$ such that $\gamma >\frac{1}{2s_2}$, and $\alpha > \max\{\frac{\gamma^2}{\bar{\mu}}, \frac{1}{2s_2}\}$. Then we have $\dot{V}^s_1<0$, which leads to the convergence of algorithm \eqref{cd}.

Next, we estimate the convergence rate of \eqref{cd} when $\Omega_i=\mathbf{R}^m, i=1,\cdots,n$.
In this case $\Pi_{\Omega_i} (x_i, -\nabla f_i(x_i) + \lambda_i) =  -\nabla f_i(x_i) + \lambda_i$,
and $\beta(x_i)n_{i}(x_i)=\mathbf{0}$. Still take $V_1^s$ in \eqref{vs1}, and then \eqref{v1dot} still holds
in this case.

Take $V^s_2= \varepsilon (\theta-\eta)^T(\theta-\eta)$, and we have
\begin{equation}
\begin{array}{lll} \label{v2d}
\dot{V}^s_2 &=&   - \varepsilon Y^T h + \varepsilon \theta^T \theta + \varepsilon \theta_2 R^T L R \eta_2 + \varepsilon \theta_2 R^T L R \delta_2 \\
          &&+   \varepsilon \eta^T [r, R]^T  h- \varepsilon \eta^T \eta- \varepsilon \eta_2 R^T L R \eta_2- \varepsilon \eta_2 R^T L R \delta_2 \\
          &\leq& - (\varepsilon \overline{\mu}-\varepsilon) \theta^T \theta +  \frac{1}{2} \varepsilon s_n^2 \theta^T \theta
           + \frac{1}{2}\varepsilon \eta_2^T\eta_2^T
          +  \frac{1}{2} \varepsilon s_n^2 \theta^T \theta \\
&&+ \frac{1}{2} \varepsilon \delta_2\delta_2^T
          +  \frac{1}{2} \varepsilon \eta^T\eta + \frac{1}{2} \varepsilon M^2 \theta^T \theta - \varepsilon \eta^T \eta- \varepsilon s_2\eta^T_2 \eta_2\\
          &&+  \frac{1}{2} \varepsilon s_n^2 \eta_2^T\eta_2 +\frac{1}{2} \varepsilon \delta^T_2\delta_2 \\
          &\leq& - \varepsilon( \overline{\mu}-1- s_n^2-\frac{1}{2} M^2 ) \theta^T\theta-\frac{1}{2}\varepsilon\eta^T\eta\\
          &&+ \varepsilon \delta_2^T\delta_2- \varepsilon( s_2- \frac{1}{2}- \frac{1}{2} s_n^2 )\eta_2^T\eta_2, \nonumber
\end{array}
\end{equation}
with $M=\max \{k_1,\cdots,k_n\}$, by using the inequality $ x^Ty \leq \frac{1}{2} ||x||_2^2+|| y||_2^2 $ and \eqref{lap} in the first step of \eqref{v2d}.

With $V^s=V^s_1+V^s_2$, it is easy to see that
\begin{equation}
\begin{array}{lll}
\dot{V}^s&=&  -(\alpha \bar{\mu} -\gamma^2 + \varepsilon( \overline{\mu}-1- s_n^2-\frac{1}{2} M^2 ) ) \theta^T \theta  \\
       &-& \frac{1}{2} \varepsilon \eta^T\eta -(\gamma s_2 -\frac{1}{2}- \varepsilon ) \delta_2^T \delta_2\\
       &-& (\alpha s_2 - \frac{1}{2}+ \varepsilon( s_2- \frac{1}{2}- \frac{1}{2} s_n^2 )) \eta_2^T\eta_2. \\
\end{array}
\end{equation}

Choose $\gamma \geq \frac{ 3\varepsilon +1}{2s_2}$ such that $\gamma s_2 -\frac{1}{2}- \varepsilon \geq \frac{1}{2}\varepsilon $.
Select $\alpha$ such that
\begin{equation}\label{par111}
\alpha \bar{\mu} -\gamma^2 + \varepsilon( \overline{\mu}-1- s_n^2-\frac{1}{2} M^2 ) \geq  \frac{1}{2} \varepsilon,
\end{equation}
 and
\begin{equation}\label{par222}
\alpha s_2 - \frac{1}{2}+ \varepsilon( s_2- \frac{1}{2}- \frac{1}{2} s_n^2 ) \geq 0.
\end{equation}
As a result, $$ \dot{V}^s \leq \frac{1}{2} \varepsilon ||p ||_2^2.$$
Then, with $
\frac{1}{2} \alpha ||p ||_2^2 \leq  V^s \leq (\frac{1}{2}(\alpha+\gamma)+\gamma + 2\varepsilon) ||p ||_2^2
$, we have
\begin{equation}\label{par12}
||p|| \leq \sqrt{ \frac{\alpha+3\gamma + 4\varepsilon}{\alpha}} || p(0)||e^{-\frac{2\varepsilon}{\alpha+3\gamma + 4\varepsilon}t},
\end{equation}
which leads to the exponential convergence of  algorithm \eqref{cd} to its equilibrium point.
\hfill $\Box$

\begin{rem}
In fact, the  exponential convergence speed  $ \frac{2\varepsilon}{(\alpha+3\gamma + 4\varepsilon)}$ in \eqref{par12} can be estimated
by solving the following optimization problem
\begin{equation}
\max_{\alpha, \gamma, \varepsilon \geq 0}    \frac{2\varepsilon}{(\alpha+3\gamma + 4\varepsilon)}
     \; \qquad  s. \; t. \;    \gamma \geq \frac{ 3\varepsilon +1}{2s_2}, \;  \eqref{par111},\; \eqref{par222}. \nonumber
\end{equation}
To get a simple estimation, we take $\gamma=\frac{3\varepsilon +1}{2s_2}$.
 Denote $\varrho_1=(s_n^2+1-2s_2)$ and $\varrho_2= s_2^2(6+4s_n^2+2M^2-4\bar{\mu}) $.  With taking
$\alpha = \max \{ \frac{1+\varepsilon\varrho_1}{2s_2}, \frac{ 9\varepsilon^2 + \varepsilon(6+\varrho_2)+1 }{4s_2^2\bar{\mu}} \} $, we have
$ \frac{2\varepsilon}{(\alpha+3\gamma + 4\varepsilon)}
 \geq  \min \{ \frac{2 s_2}{8+6s_2+s_n^2}, \frac{ 4s_2^2 \bar{\mu} }{ (3+2s_n^2+M^2+6\bar{\mu}) s_2^2+ 9\bar{\mu} s_2 +3\sqrt{6\bar{\mu}s_2+1} + 3  }  \}$.
\end{rem}

\section{Distributed economic dispatch in power grids}

In this section, the algorithm proposed in  Section 4
is applied  to the DEDP in power grids to find the optimal secure
generation allocations for power balancing  in a distributed manner.
Example \ref{exa2} is given to show that the distributed algorithm
can efficiently adapt to online network data/configuration changes,
including  {\bf load demands, generation costs/capacities, and
plugging-in/leaving-off of buses}, while Example \ref{mexa} with a
large-scale network  illustrates the {\bf scalability} of the
proposed algorithm.

Suppose that there exist control areas $\mathcal{N}=\{1,...,n\}$ with area $i$
having local generators  to supply power $P_i^g\in \mathbf{R}_{\geq 0}$
and local load demands $P_i^d \in \mathbf{R}_{\geq 0}$ to be met.
The local generation must be kept within the capacity or security bounds $ \underline{P}_i \leq P_i^g \leq \bar{P}_i, \underline{P}_i, \bar{P}_i\in \mathbf{R}_{\geq 0}$.
$f_i(P_i^g): \mathbf{R}_{\geq 0}\rightarrow \mathbf{R}$  represents the local generation cost  in control area $i$ with respect to its local generation $P_i^g$,
and it satisfies Assumption \ref{asumFun}.
Then the DEDP formulation can be derived in the form of (\ref{CRA1})
\begin{equation}\label{edp}
\begin{array}{l}\hline
{\bf Distributed \;  Economic \; Dispatch \; Problem }:           \\ \hline
\displaystyle  \min_{P_i^g,i\in \mathcal{N}}  \;\; \;f(P^g)= \sum_{i\in \mathcal{N}} f_i(P_i^g)           \\
\displaystyle  subject \; to \; \sum_{i\in \mathcal{N}} P_i^g =\sum_{i\in \mathcal{N}} P_i^d  \\
\displaystyle \qquad \qquad    \underline{P}_i \leq  P_i^g \leq  \overline{P}_i, \quad \mathcal{N}=\{1,...,n\}                             \\  \hline
\end{array}
\end{equation}

Here a multi-agent network is introduced to solve the DEDP  (\ref{edp}) motivated by recent DEDP works like \cite{ED3}, \cite{EDcc}, \cite{ED1}, \cite{ED_zam} and \cite{ED2}.
Agent $i$ is responsible to decide the generation $P^g_i$ in control area $i$
to minimize the global cost as the sum of all individuals' generation costs, while meeting the total load demands within its capacity bounds.
In addition,  each agent can  react to the changing local environment in real time,
and adapt its own behavior or preference by adjusting its local data, including $P_i^d, \underline{P}_i,  \bar{P}_i, f_i(P_i^g)$.
The  agents can also share  information
with their neighbors to facilitate the cooperations.

Then applying \eqref{cd} to \eqref{edp}, the distributed algorithm for agent $i$ is
\begin{equation}\label{cdedp}
\begin{array}{l}\hline
{\bf Distributed \;  algorithm \; for \; agent \; i}:           \\ \hline
\displaystyle \dot{P}^g_i      =[-\nabla f_i(P_i^g)+\lambda_i]^{\bar{P}_i-P_i^g}_{P_i^g-\underline{P}_i}            \\
\displaystyle \dot{\lambda}_i  =-\sum_{j\in \mathcal{N}_i}(\lambda_i-\lambda_j)-\sum_{j\in \mathcal{N}_i}(z_i-z_j)+(P^d_i-P^g_i)  \\
\displaystyle \dot{z_i}        =\sum_{j\in \mathcal{N}_i}(\lambda_i-\lambda_j)                            \\  \hline
\end{array}
\end{equation}
where $[v]_{x-c_1}^{c_2-x}=0$ if $x-c_1=0$ and  $ v \leq 0$ or $c_2-x=0$ and $v \geq 0$, otherwise $[v]_{x-c_1}^{c_2-x}=v$.

The algorithm  \eqref{cdedp}  ensures generation capacity bounds explicitly, and converges
without  the initialization procedure, which is crucially important for  the ``plug-and-play" operation in the future smart grid.

\begin{rem}
The optimization problem  \eqref{CRA1} can also be applied to  model the multi-period demand response in power systems. The objective functions describe the dis-utility of cutting loads in each area.
The resource constraint specifies the amount of loads  to be shed in the multi-periods. Particularly, the local load shedding constraints concern with the
total power demands in the multi-periods and other specifications of that area.
Then, the previous algorithms \eqref{gcd} and \eqref{cd} can be applied to solve this multi-dimensional DRAO
with general  LFCs in a distributed manner.
This issue is interesting but beyond the scope of this paper, and will be discussed elsewhere.
\end{rem}

The following two examples are presented to further show the online
adaptation property and scalability of our algorithm.
Firstly,  the standard  IEEE 118-bus system is adopted to illustrate the performance of
algorithm \eqref{cdedp}.

\begin{exa}[Optimality and adaptability]\label{exa2}
 Consider  the DEDP \eqref{edp} in the IEEE 118-bus
system with  $54$ generators. Each generator has a local quadratic generation
cost function as $f_i(P_i^g)=a_i {P_i^g}^2 +b_i P_i^g+c_i$, whose
coefficients belong to the intervals $a_i\in [0.0024,0.0679](M\$/MW^2)$,
$b_i\in [8.3391,37.6968](M\$/MW)$, and $c_i\in [6.78,74.33](M\$)$. The generation
capacity bounds of the generators are drawn from $\underline{P}_i\in
[5,150](MW)$ and $\bar{P}_i\in [150,400](MW)$, while the load of each bus
varies as $P_i^d\in [0,300](MW)$.

The corresponding agents share information on an undirected ring graph with
additional undirected edges $(1,4)$, $(15,25)$, $(25,35)$, $(35,45)$
and $(45,50)$. The simulations are performed with the differentiated
projected algorithm \eqref{cdedp}. The initial generation $P_i^g$ in
algorithm (\ref{cdedp}) is set within its local capacity bounds,
while  variables $\lambda_i$'s and $z_i$'s are set with zero initial values.

Next explains how the network data/configuration changes
at different times.
\begin{itemize}
\item {\bf Load variations:} At time 100s, 18 buses are randomly chosen
with randomly varying their loads by $-20 \sim +20\%$.
\item {\bf Generation capacity variations:} At time 200s, 18 generators
randomly vary their capacity lower bounds by  $-50 \sim +50\%$, and another 18
generators randomly vary their capacity upper bounds by $-20
\sim +20\%$.
\item {\bf Generation cost variations:}  At time 300s, 18 generators randomly
vary their $a_i$ by $0 \sim 50\%$, and another 18 generators randomly vary  their  $b_i$ by  $-50 \sim 0\%$.
\item {\bf Leaving-off of buses:} At time 400s, generator 2 and 3 disconnect from the system, and the communication edges associated to them are also removed.
\item {\bf Plugging-in of bus:}  At time 500s, generator 3 plugs in the system
with re-generated  configurations
$a_i,b_i,c_i,P_i^d,\underline{P}_i$, and $\bar{P}_i$. The undirected
edge $(3,4)$ is also added to the communication graph.
\end{itemize}

\begin{figure}
  \centering
  \includegraphics[width=3.1in]{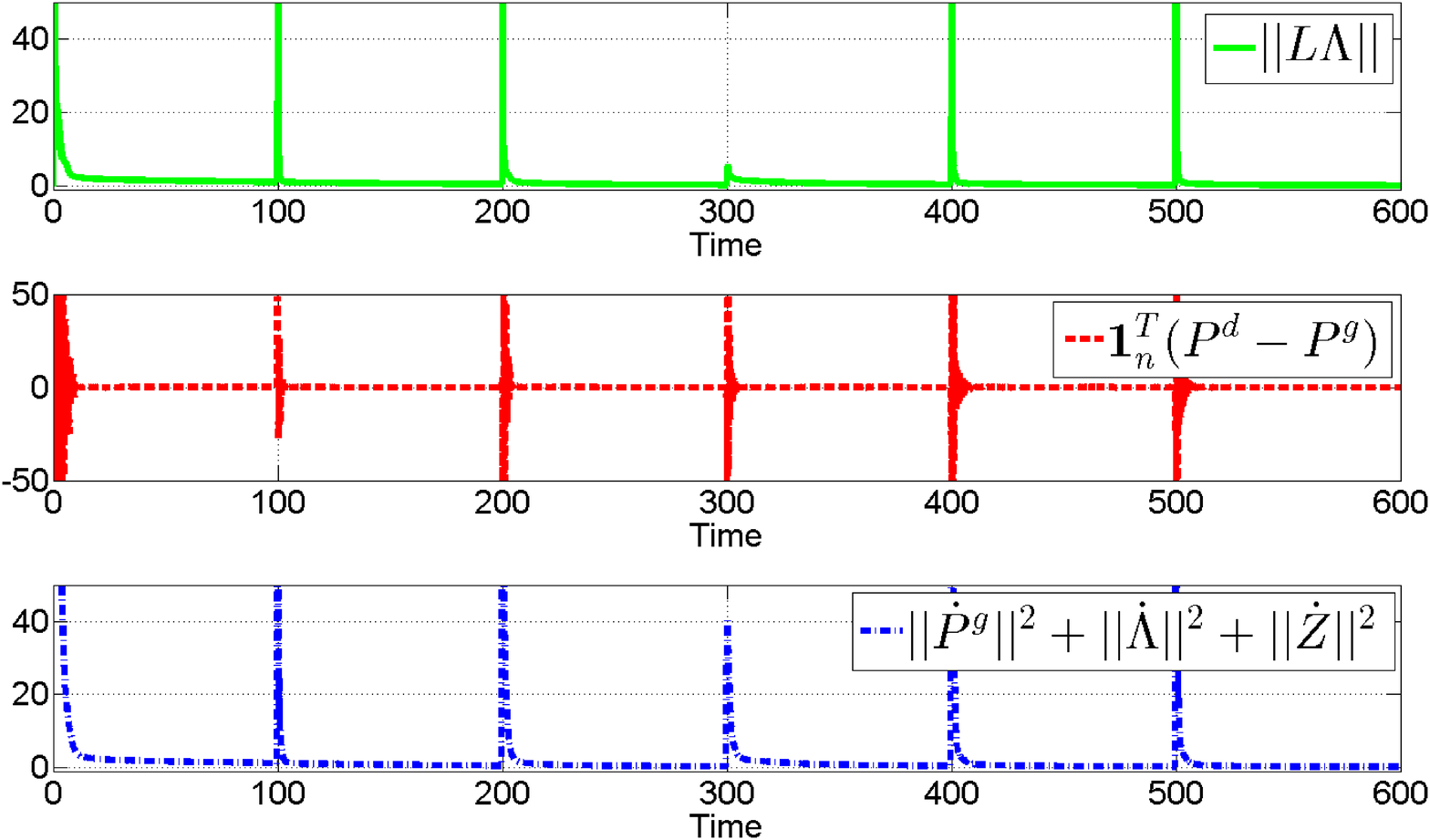}\\
  \caption{Algorithm performance indexes:
(i), $||L\Lambda||$'s always decrease to zero even with different $L$'s, implying that $\lambda's$ always reach consensus.
(ii),
The power balance constraint (i.e., resource constraint) is violated if any network configuration changes.
But the power balance gap $\mathbf{1}^T_n(P^d-  P^g)$  asymptotically decreases to zero,
even without any re-initialization coordination.
(iii), The optimality condition $||\dot{P}^g ||^2 + || \dot{\Lambda} ||^2 +  || \dot{Z}||^2 =0$
can always be satisfied asymptotically, implying the economic efficiency of the generation dispatch.}\label{fig2}
\end{figure}

When the data/configuration changes, each agent only projects its local generation onto to its local capacity bounds if necessary.
The trajectrories of  dynamics \eqref{cdedp} are derived with a first-order Euler discretization using Matlab.
The trajectories of some algorithm performance indexes are shown in Figure  \ref{fig2}.
It indicates that algorithm \eqref{cdedp} can  adaptively  find the optimal solutions to \eqref{edp} in a fully  distributed way, even without any initialization coordination procedure, when the network data or configuration  changes.
\end{exa}

\begin{exa}[Scalability]\label{mexa}
This example considers a  network of ${\bf 1000}$  control areas to achieve economic dispatch  during a normal day.
Control area $i$  has cost function  $f_i(P^g_i)= a_i {P^g_i}^2 +b_i P^g_i$ as well as generation capacity upper bound $\bar{P}_i$ and lower bound $\underline{P}_i$.
The control areas are divided equally into two groups. The first group, named as fuel group,
 is mainly supported with traditional thermal  generators,
 and has relatively higher generation costs  and larger  capacity ranges
with  the nominal values of   $a_i,b_i,\underline{P}_i,\bar{P}_i$ randomly drawn from the intervals, $[3,7](M\$/MW^2) $,\;
$[5,9](M\$/MW), $\;$ [2,6](MW), $\;$ [15,23](MW)$, respectively.
The second  group, named as renewable group, is mainly supported  with renewable energies,
and  has lower generation costs and smaller capacity ranges with the nominal values of   $a_i,b_i,\underline{P}_i,\bar{P}_i$ randomly drawn from
the intervals, $[\frac{1}{2},2](M\$/MW^2)$, \;
$[\frac{1}{2},4](M\$/MW) $, \;
$ [0,1](MW) $, \;
$[\frac{3}{2},7](MW)$, respectively.

\begin{figure}
  \centering
\includegraphics[width=3.1in]{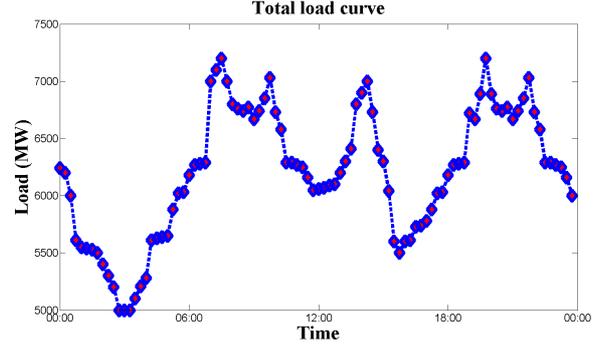}
\caption{The total load curve}\label{fig_mas1}
\end{figure}

The daily 96-point load data (15 minutes for each period) of each control
area is generated from a typical load curve for a distribution
system added with certain random perturbations. The total load curve
of the network is given in Fig \ref{fig_mas1}. In each period,
10\% of the control areas in the renewable  group  are randomly chosen to
change their generation costs and capacities  like Example
\ref{exa2} with the variations less that $\pm20\%$ of its nominal
values.

In each period, a connected graph is re-generated with random graph model  $\mathcal{G}(1000,\mathbb{P})$ as the information sharing graph of that period.
In $\mathcal{G}(1000,\mathbb{P})$,  every possible edge occurs independently with the probability of  $\mathbb{P}$.
Here,  the probability $\mathbb{P}$ in each period is  randomly drawn from the interval $[0.0015,0.005]$.

For each period, the computation time is set as $80s$. Figure \ref{fig_mas2} shows the histogram of consensus error $||L\Lambda||_2$,
 power balance gap $\sum P_i^g -\sum P_i^d$ and  optimality condition $||\dot{P}^g||_2+||\dot{\Lambda}||_2+||\dot{Z}||_2$ at time $80s$.
 It indicates that the agents  can always find the economic power  dispatch with varying loads and generation costs/capacities,  and evidently demonstrates the scalability of the proposed method.
\end{exa}

\begin{figure}
  \centering
  \includegraphics[width=3.1in]{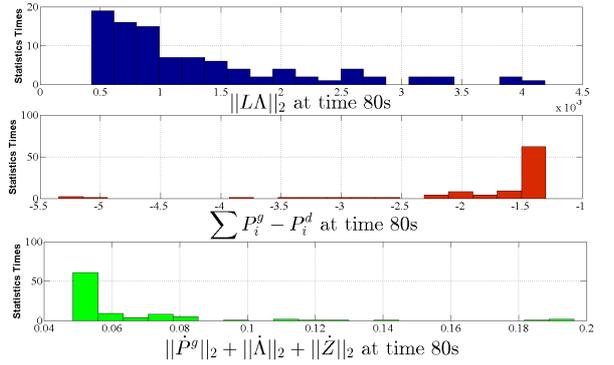}\\
  \caption{ The performance indexes at time $t=80s$ in histograms: (i), The consensus error $||L\Lambda||_2$ at $t=80s$ always decreases to
 a rather lower level due to the second-order proportional-integral consensus dynamics.
(ii),  The power balance is almost  achieved at $t=80s$ without the initialization coordination.
(iii), The optimality condition can always be satisfied at $t=80s$.
  }\label{fig_mas2}
\end{figure}

\section{Conclusions}

In this paper, a  class of  projected  continuous-time distributed  algorithms have been proposed  to solve resource allocation optimization problems with the consideration of LFCs.
The proposed algorithms  are scalable and  free of initialization coordination procedure,
and therefore, are adaptable  to working condition variations.
These salient features have important implications in the DEDP in power systems.
Firstly it allows  quite general non-box LFCs,
which is crucial when inverter-based devices are involved because their LFCs are usually  depicted in a quadratic form.
Secondly, our method is initialization free,
which may facilitate the implementation of the so-called ``plug-and-play" operation
for future smart grids in a dynamic environment.
Our application examples illustrate such implications, showing an appealing potential in the smart operation of future power grids.

We would like to note that  many challenging DRAO problems
still remain to be investigated, including the design of algorithms for the non-smooth objective functions based on differential inclusions,
the estimation of convergence rates for the proposed algorithms with general LFCs, and the development of stochastic algorithms to achieve the DRAO with noisy data observations.
Furthermore, inspired by \cite{ED1}, we also hope to combine our algorithms with physical dynamics of
power grids to derive a more comprehensive solution for the DEDP in power systems.

\end{document}